\documentclass[runningheads]{llncs}
\usepackage[T1]{fontenc}
\usepackage{graphicx}

\usepackage{color}

\usepackage{amsmath,amsfonts}
\usepackage{cleveref}
\graphicspath{{./Figures/}}
\newif\ifplots
\plotstrue
\usepackage{mathtools}
\usepackage{multirow}
\usepackage{todonotes}
\usepackage{enumitem}
\usepackage{pgfplots}
\pgfplotsset{compat=newest}
\usetikzlibrary{pgfplots.groupplots}

\usepackage{booktabs}
\usepackage{makecell}

\usepackage{tikz}
\newcommand{\R}{\mathbb{R}}

\DeclareMathOperator{\diag}{diag}

\begin{document}
\title{Application of Deep Kernel Models for Certified and Adaptive RB-ML-ROM Surrogate Modeling}
\titlerunning{Deep Kernel Models for RB-ML-ROM Modeling}
\author{Tizian Wenzel\inst{1}\thanks{Corresponding author, \email{tizian.wenzel@mathematik.uni-stuttgart.de}} \and
Bernard Haasdonk\inst{1} \and
Hendrik Kleikamp\inst{2} \and
Mario Ohlberger\inst{2} \and
Felix Schindler\inst{2}}
\authorrunning{T. Wenzel et al.}
\institute{Institute for Applied Analysis and Numerical Simulation,\\ University of Stuttgart, Germany\\ \email{\{tizian.wenzel,bernard.haasdonk\}@mathematik.uni-stuttgart.de} \and Institute for Analysis and Numerics, Mathematics Münster,\\ University of Münster, Germany\\ \email{\{hendrik.kleikamp,mario.ohlberger,felix.schindler\}@uni-muenster.de}}
\maketitle              %
\begin{abstract}
In the framework of reduced basis methods, we recently introduced a new certified hierarchical and adaptive surrogate model, which can be used for efficient approximation of input-output maps that are governed by parametrized partial differential equations.
This adaptive approach combines a full order model, a reduced order model and a machine-learning model. 
In this contribution, we extend the approach by leveraging novel kernel models for the machine learning part,
especially structured deep kernel networks as well as two layered kernel models.
We demonstrate the usability of those enhanced kernel models for the RB-ML-ROM surrogate modeling chain and highlight their benefits in numerical experiments.

\end{abstract}
\section{Introduction} \label{sec:introduction}

Model Order Reduction (MOR) is an indispensable technology in the context of multi-query, real-time or slim computing context,
as full oder models (FOMs) of real  physical processes may be unaffordable.
Many techniques for MOR exist, cf.\ \cite{benner2020model}.
A large class of MOR methods are reduced basis (RB) models, which are projection-based and the projection spaces are typically constructed from samples of the FOM and allow rigorous a-posteriori error estimation by residual analysis~\cite{haasdonk2017reduced}.
As a drawback, the time-integration requires time-marching which may be a bottleneck in case of a fine time resolution.
In recent years, also machine learning (ML) based techniques have entered the field, which enable rapid prediction without expensive time-marching \cite{wang2019nonintrusive}.
Such data-based approaches, however, mostly lack any error guarantees.
Recently, we presented an adaptive and hierarchical RB-ML-ROM framework, which combines the benefits of those two approaches, i.e.~rapid prediction of ML-models including certification from RB methods \cite{haasdonk2022certified}.
In the current presentation, we demonstrate the versatility of the RB-ML-ROM framework, by including efficient multilayer kernel models for the ML part, namely Structured Deep Kernel Networks (SDKN) and two-layered greedy kernel models (2L-VKOGA).

The paper is organized as follows: In \Cref{subsec:rom} and \Cref{subsec:ml} we recall some more background information on reduced order modeling and machine learning. 
\Cref{sec:deep_kernels} elaborates on novel deep kernel models and \Cref{sec:certified} reviews the recently introduced certified and adaptive RB-ML-MOR model.
Subsequently, \Cref{sec:num_appl} provides the numerical experiments on the combination of these approaches and finally \Cref{sec:conclusion} summarizes and gives an outlook.
\vspace{-.3cm}

\section{Reduced order modeling} \label{subsec:rom}
\vspace{-.3cm}

In this paper, we consider parametrized parabolic PDEs on a bounded Lipschitz-domain~$\Omega\subset\R^m$ over the time interval~$[0,T]$ with~$T>0$.
Further, let~$\mathcal{P}\subset\R^p$ be the set of parameters.
We denote by~$V_h$ a Hilbert space associated to a finite element discretization of~$\Omega$, such that~$V_h\subset H^1(\Omega)\subset L^2(\Omega)$.
For a parameter~$\mu\in\mathcal{P}$, we consider the semi-discrete problem of finding~$u_h(\mu) \in L^2(0,T; V_h)$ with $\partial_t u_h(\mu) \in L^2(0,T; V_h')$ such that
\begin{align}
	\langle\partial_t u_h(\mu), v_h\rangle + a\big(u_h(\mu), v_h; \mu\big) &= l(v_h; \mu) &&\forall \ v_h \in V_h \ \text{in } [0,T],\label{equ:fom-problem}\\
	u_h(0; \mu) &= u_0(\mu), &&\notag
\end{align}
where~$u_0(\mu)\in V_h$ is the initial condition.
In addition, we are particularly interested in output trajectories, i.e.~we consider an output functional~$f_h\colon\mathcal{P}\to L^2(0,T)$ defined as~$f_h(t;\mu)\coloneqq s(u_h(t;\mu))$ for some~$s\in V_h'$.
For every~$\mu\in\mathcal{P}$, $l(\,\cdot\,;\mu)\in V_h'$ is assumed to be a linear and continuous functional and~$a(\,\cdot\,,\,\cdot\,; \mu)\colon V_h \times V_h \to \R$ denotes a continuous and coercive bilinear form.
\par
Since solving the problem above for many different values of the parameter is costly, we are looking for a reduced order model that is faster to evaluate.
To this end, the reduced basis (RB) method projects the problem~\eqref{equ:fom-problem} onto a well-chosen reduced space~$V_\text{RB}\subset V_h$ with~$\dim V_\text{RB}\eqqcolon N_\text{RB}\ll N_h\coloneqq\dim V_h$ (for a general introduction to RB methods, see~\cite{haasdonk2017reduced}).
One thus searches for an approximate solution~$u_\text{RB}(t;\mu)\in V_\text{RB}$ of~$u_h(t;\mu)\in V_h$ by solving the same semi-discrete problem as in~\eqref{equ:fom-problem} with test functions from the space~$V_\text{RB}$, i.e.
\begin{align}
	\langle\partial_t u_\text{RB}(\mu), v_\text{RB}\rangle\hspace*{-2pt} +\hspace*{-1pt} a\big(u_\text{RB}(\mu), v_\text{RB}; \mu\big) &= l(v_\text{RB}; \mu) &&\forall \ v_\text{RB} \in V_\text{RB}  \  \text{in } [0,T] ,\label{equ:rom-problem}\\
	\langle u_\text{RB}(0;\mu),v_\text{RB}\rangle &= \langle u_0(\mu),v_\text{RB}\rangle \hspace*{-5pt}&&\forall \ v_\text{RB} \in V_\text{RB}.\notag
\end{align}
To obtain an efficient reduced order model, we pose the standard assumption that the bilinear form~$a$ and the right hand side~$l$ are parameter separable.
In that case, many of the required quantities can be precomputed during the so-called \emph{offline}-phase.
In the subsequent \emph{online}-phase, the reduced solution for a new parameter can be determined independently of the dimension~$N_h$ of~$V_h$.
\par
To assess the error between~$u_h$ and~$u_\text{RB}$, one can consider the residual by inserting~$u_\text{RB}$ into Eq.~\eqref{equ:fom-problem}.
The norm of the residual can be computed without actually solving~\eqref{equ:fom-problem}.
Further, the corresponding error estimator~$E_\text{RB}(u;\mu)$ for arbitrary~$u\in V_\text{RB}$, introduced in~\cite{grepl2005aposteriori}, can be evaluated efficiently during the online phase if~$a$ and~$l$ fulfill the parameter separability assumption.
Since the functional~$s$ is assumed to be linear and continuous, the error estimator for the state~$u_\text{RB}$ also results in an estimator for the error in the output functional, i.e.
\[
\lVert f_h(\mu)-s(u_\text{RB}(\mu))\rVert_{L^2(0,T)}\leq\lVert s(\,\cdot\,;\mu)\rVert_{V_h'}E_\text{RB}(u_\text{RB}(\mu);\mu).
\]
We emphasize at this point that the error estimator computes an estimate based on the parameter~$\mu$ and coefficients with respect to the reduced basis.
However, these coefficients do not have to be the result of solving the problem in~\eqref{equ:rom-problem}.
\vspace{-.3cm}

\section{Machine learning for regression} \label{subsec:ml}
\vspace{-.3cm}

In the case of regression, ML procedures typically assume that input data $X_N \coloneqq \{x_1, \dots, x_N \} \subset \R^d$ and corresponding output data $Y_N \coloneqq \{ y_1, \dots, y_N \} \subset \R^b$ is given.
The task is to ''learn`` a function which approximates the mapping $x_i \mapsto y_i$ and generalizes properly to unseen novel inputs $x$.
One mathematically well established approach for this task is to employ kernel methods.
These kind of methods can be understood as applying linear regression to a nonlinear transformation of the data, which is implicitly given by the kernel~$k$, which is in our case a symmetric function $k\colon \R^d \times \R^d \rightarrow \R$. 
A well known representer theorem then states that $f^*$ can be expanded as
\begin{align}
	\label{eq:kernel_model}
	s_{X_N} = \sum_{j=1}^N \alpha_j k(\cdot, x_j).
\end{align}
In order to obtain sparse models, greedy methods can be employed to select a suitable subset of centers $X_n \subset X_N$ such that the kernel expansion from Eq.~\eqref{eq:kernel_model} is smaller.
An algorithm for this, which will be leveraged later on, is the vectorial kernel orthogonal greedy algorithm (VKOGA) \cite{santin2021kernel}. \\
Another lass of ML methods is deep learning, where probably neural networks are the most prominent example.
Due to their deep structure, feature learning is possible which can be beneficial especially for medium to high dimensional datasets. \\
In order to make the best out of both approaches, the idea is to combine deep learning methods such as neural networks with the benefits of ''shallow`` kernel methods.
In the following \Cref{sec:deep_kernels}, we review two recently introduced algorithms in this direction:
SDKNs \cite{wenzel2021universality} as well as 2L-VKOGA \cite{wenzel2023data}.
\vspace{-.3cm}

\section{Deep kernel models} \label{sec:deep_kernels}
\vspace{-.3cm}
Both, the SDKN and the 2L-VKOGA approach, are based on a deep kernel representer theorem \cite{bohn2019representer},
which allows to put kernel methods into a multilayer setup in a mathematical principled way: \\
\textbf{SDKNs} \cite{wenzel2021universality} make use of special classes of kernels, which are used alternately in a structured way to derive a powerful setup:
First, matrix-valued linear kernels are used, i.e.\ $k_\text{lin}(x, z) = \langle x, z \rangle_{\R^{d_0}} \cdot I_{d_1}$ with $I_{d_1} \in \R^{d_1 \times d_1}$ being the identity matrix.
Second, single-dimensional kernels are employed that act on single components $x^{(i)}, z^{(i)}$ separately as $k_s(x, z) = \diag(k(x^{(1)}, z^{(1)}), \dots, k(x^{(d)}, z^{(d)}))$. 
Hereby, the one dimensional kernel $k\colon \R \times \R \rightarrow \R$ is e.g.\ the Gaussian kernel.
The SDKN model is set up using the linear kernels for every odd layer and single-dimensional kernels for every even layer.
Thus, the odd layers essentially allow for linear combinations of the inputs, while the even layers can be compared to activation functions of neural networks.
However, due to being kernel models that can be optimized instead of fixed activation functions, they allow for more flexibility and thus enable a potentially faster optimization. %
As analyzed in \cite{wenzel2021universality}, these SDKNs enjoy universal approximation properties in various limit cases.\\
The \textbf{2L-VKOGA} \cite{wenzel2023data} combines the use of an optimized two-layered kernel with greedy kernel algorithms, where we leverage the VKOGA algorithm \cite{santin2021kernel}: %
Radial basis function kernels are a popular class of kernels which are given as $k(x, z) = \Phi(\Vert x - z \Vert)$
for some radial basis function $\Phi\colon \R^d \rightarrow \R$,
e.g.\ the Gaussian kernel $k(x, z) = \exp(-\Vert x - z \Vert^2)$.
These RBF kernels usually make use of a shape parameter $\varepsilon > 0$ via $k_\varepsilon(x, z) = \Phi(\varepsilon \Vert x - z \Vert)$,
which can be used to adopt the shape of the kernel to the dataset. %
The two-layered kernel approach generalizes this scale parameter by first allowing for a $d \times d$ matrix $A$ instead of only a scalar valued shape parameter~$\varepsilon$, i.e.\ $k_A(x, z) \coloneqq \Phi(\Vert A(x-z) \Vert)$,
and secondly providing an ML-based gradient descent optimization, which allows to optimize the kernel (i.e.\ the matrix $A$) efficiently.
In particular, this setup can be framed within the deep kernel representer theorem \cite{bohn2019representer} and thus be seen as a two-layered kernel. %
The optimized kernel $k_A$ can subsequently be used with standard kernel algorithms.
Here we use it in conjunction with greedy algorithms, and in particular with the VKOGA algorithm:
In view of Eq.~\eqref{eq:kernel_model}, the VKOGA algorithm selects a suitable subset $X_n \subset X_N$ of the dataset, such that the resulting kernel model consists only of $n \ll N$ summands and is thus sparse. %
For our current application, we extend the 2L-VKOGA approach that was only introduced for scalar valued target data \cite{wenzel2023data} to vectorial data, i.e. $Y_N \subset \R^b$.
Therefore, we consider for the kernel optimization the $\ell^2$-norm over the errors from the single dimensions.
As the VKOGA algorithm is intrinsically already vector-valued, both can be combined.

\vspace{-.3cm}
\section{Certified and adaptive RB-ML-ROM} \label{sec:certified}
\vspace{-.3cm}

We summarize the adaptive model hierarchy as specified in~\cite{haasdonk2022certified}.
Instead of describing the model hierarchy in an abstract way as done in~\cite{haasdonk2022certified}, we make direct use of the RB model from \Cref{subsec:rom} and the deep kernel methods from \Cref{sec:deep_kernels}.
\par
The RB-ML-ROM is built as a hierarchy consisting of three models, a full-order model (FOM), a reduced basis reduced order model (RB-ROM), and a machine learning reduced order model (ML-ROM).
As discussed, the error estimator of the RB-ROM (see \Cref{subsec:rom}) can be used to certify solutions from the RB-ROM as well as from the ML-ROM.
The guiding assumption is that the ML-ROM is the fastest model in the hierarchy.
A query to the RB-ML-ROM now proceeds as follows:
For a new parameter the ML-ROM is queried and its prediction accepted if its error estimation is sufficiently accurate, i.e.~the error is below a desired tolerance~$\varepsilon>0$.
Otherwise, the slightly more expensive RB-ROM is queried and its error is again validated using the error estimator.
Only if the RB-ROM is also inaccurate, one has to fall back to the costly FOM, which is assumed to be accurate enough for every parameter.
\par
Not only the prediction of the hierarchical model is adaptive, but the submodels themselves are adaptively updated as follows: When querying the adaptive RB-ML-ROM for the first time for a certain parameter~$\mu\in\mathcal{P}$, the RB-ROM and the ML-ROM are initialized as empty models, i.e.~they can only return zero as the solution.
Consequently, most likely, the FOM is called for the parameter~$\mu$ and the FOM output trajectory~$f_h(\mu)$ is returned.
As intermediate quantity the FOM solution trajectory~$u_h(\mu)$ is computed, which will be used to set up the RB-ROM.
In our concrete setting, the RB-ROM is constructed using the hierarchical approximate proper orthogonal decomposition (HaPOD)~\cite{himpe2018hierarchical} applied to data obtained from querying the FOM.
Whenever solving the FOM is necessary (due to the RB-ROM and ML-ROM not being accurately enough for the current parameter), a new FOM trajectory becomes available and the RB-ROM can be improved using this new information.
Finally, the ML-ROM is built upon the updated RB-ROM and is based on the same reduced basis and error estimator as the RB-ROM.
More precisely, using an ML algorithm such as the deep kernel models described in \Cref{sec:deep_kernels}, we learn the coefficients of the RB solution with respect to the reduced basis as a function of the parameter, see also~\cite{wang2019nonintrusive}.
The required training data is generated as a byproduct when the RB-ROM is queried.
However, only those RB solutions that fulfill the accuracy requirements (evaluated by the error estimator of the RB-ROM) are added to the training set.
Recall that the error estimator of the RB-ROM takes as input a parameter~$\mu\in\mathcal{P}$ and arbitrary reduced coefficients.
In particular, these coefficients do not have to stem from the RB-ROM, but can also be the result of an ML prediction.
This fact enables the direct application of the error estimator to certify the ML-ROM independent of the actual ML model used.
\par
For the ML-ROM, two different approaches are possible which differ in the way they treat the time component.
The \emph{random-access-in-time} method passes the time instance along with the parameter as an input to the machine learning model.
To obtain the whole solution trajectory, the machine learning model is therefore queried for all different time instances.
In contrast, the \emph{time-vectorized} approach returns the whole time sequence for a given parameter, i.e.~the output contains the reduced coefficients for all time steps at once.
\vspace{-.3cm}

\section{Numerical application} \label{sec:num_appl}

\begin{figure}[hbt!]
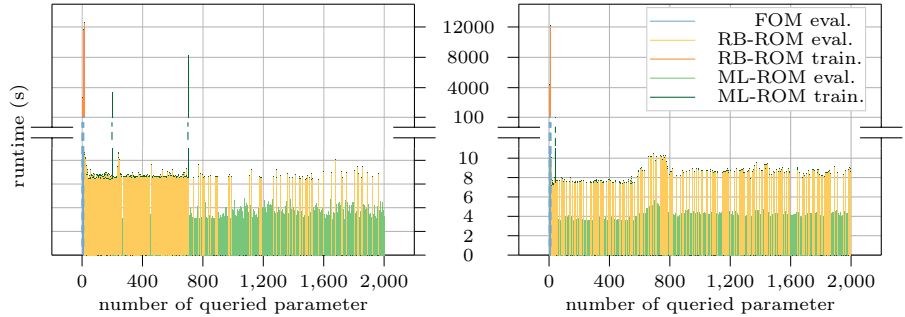

	\centering
	\begin{tikzpicture}
	\definecolor{coral25314160}{RGB}{253,141,60}
	\definecolor{cornflowerblue116169207}{RGB}{116,169,207}
	\definecolor{darkgray176}{RGB}{176,176,176}
	\definecolor{darkseagreen120198121}{RGB}{120,198,121}
	\definecolor{forestgreen010455}{RGB}{0,104,55}
	\definecolor{lightgray204}{RGB}{204,204,204}
	\definecolor{sandybrown25420492}{RGB}{254,204,92}
	
	\begin{groupplot}[
			group style={
				group size=1 by 2,
				vertical sep=0pt
			},
			xmin=0,
			xmax=2000,
			enlarge x limits=true,
			tick align=outside,
			tick pos=left,
			xmajorgrids,
			ymajorgrids,
			xtick style={color=black},
			y grid style={very thin,darkgray176},
			scaled x ticks=false,
			scaled y ticks=false,
			ytick pos=right,
			yticklabels={,,},
			yticklabel style={font=\scriptsize},
			yticklabel pos=left,
			ytick style={color=black},
			xticklabel style={font=\scriptsize},
		]
		
		\nextgroupplot[
			xtick=\empty,
			axis x line*=top,
			axis y discontinuity=parallel,
			xtick={0, 400, 800, 1200, 1600, 2000},
			xticklabels={},
			xtick style={draw=none},
			ymin=-4000,
			ymax=15000,
			ytick={100, 4000, 8000, 12000},
			height=3.5cm,
			width=6.4cm
		]
		
		\begin{scope}
			\clip (axis cs:0,100) rectangle (axis cs: 2000,14000);
			
			\input{Figures/data_SDKN_top}
		\end{scope}
		
		\nextgroupplot[
			xlabel={number of queried parameter},
			xlabel style={font=\scriptsize, yshift=4pt},
			ylabel={runtime (s)},
			ylabel style={font=\scriptsize, xshift=25pt, yshift=-5pt},
			ytick={1, 3, 5, 7, 9, 11},
			ymin=1,
			ymax=10,
			xtick={0, 400, 800, 1200, 1600, 2000},
			axis x line*=bottom,
			height=3cm,
			width=6.4cm
		]
		
		\input{Figures/data_SDKN_bottom}
	\end{groupplot}

	\draw[densely dashed, cornflowerblue116169207] (SDKN1bottom) -- (SDKN1top);
	\draw[densely dashed, cornflowerblue116169207] (SDKN2bottom) -- (SDKN2top);
	\draw[densely dashed, cornflowerblue116169207] (SDKN3bottom) -- (SDKN3top);
	\draw[densely dashed, cornflowerblue116169207] (SDKN4bottom) -- (SDKN4top);
	\draw[densely dashed, cornflowerblue116169207] (SDKN5bottom) -- (SDKN5top);
	\draw[densely dashed, cornflowerblue116169207] (SDKN6bottom) -- (SDKN6top);
	\draw[densely dashed, cornflowerblue116169207] (SDKN7bottom) -- (SDKN7top);
	\draw[densely dashed, cornflowerblue116169207] (SDKN8bottom) -- (SDKN8top);
	\draw[densely dashed, cornflowerblue116169207] (SDKN9bottom) -- (SDKN9top);
	\draw[densely dashed, cornflowerblue116169207] (SDKN10bottom) -- (SDKN10top);
	\draw[densely dashed, cornflowerblue116169207] (SDKN11bottom) -- (SDKN11top);
	\draw[densely dashed, cornflowerblue116169207] (SDKN12bottom) -- (SDKN12top);
	
	\draw[densely dashed, forestgreen010455] (SDKN13bottom) -- (SDKN13top);
	\draw[densely dashed, forestgreen010455] (SDKN14bottom) -- (SDKN14top);
\end{tikzpicture}\hspace*{-0.3cm}
	\begin{tikzpicture}
	\definecolor{coral25314160}{RGB}{253,141,60}
	\definecolor{cornflowerblue116169207}{RGB}{116,169,207}
	\definecolor{darkgray176}{RGB}{176,176,176}
	\definecolor{darkseagreen120198121}{RGB}{120,198,121}
	\definecolor{forestgreen010455}{RGB}{0,104,55}
	\definecolor{lightgray204}{RGB}{204,204,204}
	\definecolor{sandybrown25420492}{RGB}{254,204,92}

	\begin{groupplot}[
			group style={
				group size=1 by 2,
				vertical sep=0pt
			},
			xmin=0,
			xmax=2000,
			enlarge x limits=true,
			tick align=outside,
			tick pos=left,
			xmajorgrids,
			ymajorgrids,
			xtick style={color=black},
			x grid style={very thin, darkgray176},
			y grid style={very thin, darkgray176},
			scaled x ticks=false,
			scaled y ticks=false,
			ytick pos=left,
			yticklabel style={font=\scriptsize, text width=width("$12000$"), align=center},
			yticklabel pos=left,
			ytick style={color=black},
			xticklabel style={font=\scriptsize},
			legend cell align={left},
			legend style={fill opacity=1, draw opacity=1, text opacity=1, draw=lightgray204, inner xsep=1pt, inner ysep=0pt, font=\scriptsize, row sep=-3pt, xshift=-3pt},
		]
		
		\nextgroupplot[
			xtick=\empty,
			axis x line*=top,
			axis y discontinuity=parallel,
			xticklabels={},
			xtick style={draw=none},
			ymin=-4000,
			ymax=15000,
			ytick={100, 4000, 8000, 12000},
			yticklabels={$100$, $4000$, $8000$, $12000$},
			xtick={0, 400, 800, 1200, 1600, 2000},
			height=3.5cm,
			width=6.4cm
		]
		
		\begin{scope}
			\clip (axis cs:0,100) rectangle (axis cs: 2000,14000);
			
			\input{Figures/data_VKOGA_top}
		\end{scope}

		\addlegendimage{const plot, color=cornflowerblue116169207}
		\addlegendentry{\hphantom{RB-}FOM eval.}
		\addlegendimage{const plot, color=sandybrown25420492}
		\addlegendentry{RB-ROM eval.}
		\addlegendimage{const plot, color=coral25314160}
		\addlegendentry{RB-ROM train.\hspace*{-2.5pt}}
		\addlegendimage{const plot, color=darkseagreen120198121}
		\addlegendentry{ML-ROM eval.}
		\addlegendimage{const plot, color=forestgreen010455}
		\addlegendentry{ML-ROM train.\hspace*{-2.5pt}}

		\nextgroupplot[
			xlabel={number of queried parameter},
			xlabel style={font=\scriptsize, yshift=4pt},
			xtick={0, 400, 800, 1200, 1600, 2000},
			ytick={1, 3, 5, 7, 9, 11},
			yticklabels={$0$, $2$, $4$, $6$, $8$, $10$},
			ymin=1,
			ymax=12,
			axis x line*=bottom,
			height=3cm,
			width=6.4cm
		]

		\input{Figures/data_VKOGA_bottom}
	\end{groupplot}

	\draw[densely dashed, cornflowerblue116169207] (VKOGA1bottom) -- (VKOGA1top);
	\draw[densely dashed, cornflowerblue116169207] (VKOGA2bottom) -- (VKOGA2top);
	\draw[densely dashed, cornflowerblue116169207] (VKOGA3bottom) -- (VKOGA3top);
	\draw[densely dashed, cornflowerblue116169207] (VKOGA4bottom) -- (VKOGA4top);
	\draw[densely dashed, cornflowerblue116169207] (VKOGA5bottom) -- (VKOGA5top);
	\draw[densely dashed, cornflowerblue116169207] (VKOGA6bottom) -- (VKOGA6top);
	\draw[densely dashed, cornflowerblue116169207] (VKOGA7bottom) -- (VKOGA7top);
	\draw[densely dashed, cornflowerblue116169207] (VKOGA8bottom) -- (VKOGA8top);
	\draw[densely dashed, cornflowerblue116169207] (VKOGA9bottom) -- (VKOGA9top);
	\draw[densely dashed, cornflowerblue116169207] (VKOGA10bottom) -- (VKOGA10top);
	\draw[densely dashed, cornflowerblue116169207] (VKOGA11bottom) -- (VKOGA11top);
	\draw[densely dashed, cornflowerblue116169207] (VKOGA12bottom) -- (VKOGA12top);
	
	\draw[densely dashed, forestgreen010455] (VKOGA13bottom) -- (VKOGA13top);
\end{tikzpicture}
	\caption{
		Visualization of the required time (vertical axis) to build and evaluate the adaptive model $M_\text{adapt}$ for each $\mu$ chosen by the random sampling (horizontal axis).
		Left: SDKN as ML-ROM. Right: 2L-VKOGA as ML-ROM.
	}
	\label{fig:adaptive_model_results}
\end{figure}

\vspace{-.3cm}
As an application for the deep kernel methods in the framework of the adaptive RB-ML-ROM model, we leverage the Monte Carlo example introduced in \cite[Section 4.2]{haasdonk2022certified}.
Likewise, we can also use the NN used there as a baseline method for comparison.
The quantity of interest is the average (in time) temperature~$\overline{f}\colon \mathcal{P} \rightarrow \R$ within a given room of a building, represented by the domain $\Omega$.
The temperature can be computed by solving the heat equation and is depending on 28 diffusion parameters for walls, doors and heaters.
In order to determine the average and variance of this output quantity over a given probability density~$\rho$ on the 28 dimensional parameter space, we rely on a Monte Carlo estimation
\begin{align*}
	\text{E}[\overline{f}] \eqqcolon \int_\mathcal{P} \rho(\mu)\overline{f}(\mu) ~ \mathrm{d}\mu \approx \frac{1}{N_\text{mc}} \sum_{\mu \in \mathcal{P}_\text{mc}} \overline{f}(\mu),
\end{align*}
and similarly for the variance. \\
For the numerical experiment, the same setup as in \cite{haasdonk2022certified} was used, i.e.\ the full order model consists of 321206 DoFs and 999 timesteps and is thus quite expensive.
Therefore, the use of the adaptive model introduced in \Cref{sec:certified} is desirable, which is applied with a tolerance of $\varepsilon = 5 \cdot 10^{-2}$.
For the machine learning models, the SDKN followed the setup of the NN, i.e.\ four layers with each an inner dimension of 128 while using the Gaussian kernel for its single-dimensional kernel layers.
The 2L-VKOGA used a quadratic Matérn kernel and 10 epochs with a batch size of 32 and an initial learning rate of $5 \cdot 10^{-3}$ for the Adam optimizer for the optimization of the two-layered kernel. 
For the greedy algorithm, the so called $f$-greedy selection criterion was used up to an expansion size of 500 centers.
The NN and SDKN were (re)trained as soon as the data of 200 new parameter values were evaluated, the VKOGA and 2L-VKOGA models were (re)trained as soon as $40$ new data samples were available.
As soon as the performance of the ML-ROM models was sufficiently good (use of ML model at least in 60\% of the queried parameter values), no more retrainings were performed.
The results of the runs of the adaptive model are collected in \Cref{tab:overview_results} as well as \Cref{fig:adaptive_model_results}:
\Cref{tab:overview_results} summarizes the key performance numbers for the different machine learning models: 
In the column ''Training ML-ROM``, the training set size $n_\mu$ and the training duration $t$ for the ML-ROMs are listed:
It can be seen that the VKOGA models are particularly fast, which is due to their smaller training batch sizes of $40$ instead of $200$, such that they are trained already earlier.
Of particular interest is the column ''ML-ROM ratio``, which shows the ratio of ML-ROM vs.\ ROM evaluations:
For the baseline NN, this ratio increases from $0\%$ (before training) to up to $87\%$ after $3$ training stages.  
The shallow VKOGA model does not seem to be accurate at all at any time, its ML-ROM ratio never exceeds 25\% -- and is thus stopped after $2000$ queried parameter evaluations.
For the SDKN, the ML-ROM ratio is already $60.3\%$ after the first training, 
thus surpassing the NN performance after two trainings. 
Its final ratio of $94.7\%$ is large enough such that no further retraining is required anymore.
The deep 2L-VKOGA model achieves a performance of $86.7\%$ already after the very first training, which included only $40$ evaluated parameter values.
The reason for this superb performance (especially in contrast to the shallow VKOGA model) can be seen in the optimized kernel.
A singular value analysis of the matrix~$A$ of the kernel $k_A$ reveals that (quite independent of the specific run,
i.e.\ random choice of the parameter $\mu$) there is one large singular value of $2.5$, while the second largest one is already only $0.35$ and the remaining ones smaller or equal than $0.35$.
This furthermore highlights the amenability of the 2L-VKOGA for interpretability of the ML model. 
The last column ''Prediction models`` shows the number of FOM, ROM and ML-ROM evaluations (which is related to the ratios from the previous column) and the mean prediction times of these models. 
The reason for this is that the size of the ML models (in terms of the number of parameters) is approximately the same, namely approximately 80000 parameters for the NN and SDKN and approximately 100000 parameters for the VKOGA models. 
In terms of total runtime the NN took 48:57h, the SDKN 48:23h, the 2L-VKOGA 46:58h and the VKOGA was stopped after 54:01h (for only 2000 parameters). \\
The source code to carry out the experiments and generate the figures can be found online.\footnote{\url{https://gitlab.mathematik.uni-stuttgart.de/wenzeltn/paper-2023-deep-kernel-for-rb-ml-rom}}%

\vspace{-.7cm}

\begin{table}
	\centering
	\begin{tabular}{ c || c c | c || c c c }
		\toprule
		\multirow{3}{*}{ML-Model} & \multicolumn{2}{c|}{Trainings of ML-ROM} & \hphantom{i}\multirow{3}{*}{\makecell{ML-ROM\\ratio}}\hphantom{i} & \multicolumn{3}{c}{\hphantom{i}Num.\ of model evaluations\hphantom{i}} \\
		\cline{2-3}\cline{5-7}
		& \hphantom{i}\makecell{Param.\ at which\\training started}\hphantom{i} & \hphantom{i}\makecell{Training\\time (s)}\hphantom{i} & & Model & \makecell{Num.\\evals.} & \makecell{Average\\time (s)} \\
		\midrule \midrule
		\multirow{4}{*}{NN} & & & 0.0\% & FOM & 12 & $234.9$ \\
		& $200$ & $536.3\hphantom{0}$ & 36.8\% & RB-ROM & 876 & $6.7$  \\
		& $516$ & $3919.1$ & 49.1\% & ML-ROM & 4112 & $3.2$ \\
		& $908$ & $6515.1$ & 92.9\% & & & \\ \hline
		\multirow{4}{*}{VKOGA} & & & 0.0\% & FOM & 12 & $235.9$ \\
		& $40$ & $97.8$ & $\leq 25\%$ & RB-ROM & 1789 & $8.6$  \\
		& $\cdots$ & & $\cdots$ & ML-ROM & 199 & $4.2$ \\
		& $1999$ & $3572.3$ & $\leq 25\%$ & & & \\
		\midrule \midrule
		\multirow{3}{*}{SDKN} & & & 0.0\% & FOM & 12 & $233.0$ \\
		& $200$ & $3350.2$ & 60.4\% & RB-ROM & 613 & $6.7$  \\
		& $703$ & $8272.9$ & 94.7\% & ML-ROM & 4375 & $3.3$ \\ \hline
		\hphantom{i}\multirow{3}{*}{2L-VKOGA}\hphantom{i} & & & 0.0\% & FOM & 12 & $235.9$ \\
		& $40$ & $110.3$ & 86.7\% & RB-ROM & 685 & $8.7$  \\
		& & & & ML-ROM & 4303 & $4.1$ \\
		\bottomrule
	\end{tabular}
	\caption{Overview of the results using the adaptive model for different kind of machine learning models:
		Top two boxes: Baseline NN as well as a shallow kernel VKOGA model (\Cref{subsec:ml}).
		Bottom two boxes: The deep kernel approaches SDKN and 2L-VKOGA (\Cref{sec:deep_kernels}).}
	\label{tab:overview_results}
\end{table}

\vspace{-.3cm}
\section{Conclusion} \label{sec:conclusion}
\vspace{-.3cm}

In this paper we reviewed the certified RB-ML-ROM model and the two deep kernel models SDKN and 2L-VKOGA, whereby we extended the second one for the use with vectorial target data.
These deep kernel models where then used for the machine learning part within the RB-ML-ROM model, where they outperformed the standard NN baseline as well as shallow kernel models. \\

\vspace{-.3cm}
\subsubsection{Acknowledgements}
Funded by BMBF under the contracts 05M20PMA and 05M20VSA. 
Funded by the Deutsche Forschungsgemeinschaft (DFG, German Research Foundation) under Germany’s Excellence Strategy EXC 2044 --390685587, 
Mathematics Münster: Dynamics--Geometry--Structure, and EXC 2075 --390740016.
We acknowledge the support by the Stuttgart Center for Simulation Science (SimTech).

\bibliographystyle{splncs04}
\bibliography{references}

\begin{thebibliography}{10}
\providecommand{\url}[1]{\texttt{#1}}
\providecommand{\urlprefix}{URL }
\providecommand{\doi}[1]{https://doi.org/#1}

\bibitem{benner2020model}
Benner, P., Grivet-Talocia, S., Quarteroni, A., Rozza, G., Schilders, W.,
  Silveira, L. (eds.): Model Order Reduction. Volume 1. De~Gruyter, Berlin
  (2020)

\bibitem{bohn2019representer}
Bohn, B., Rieger, C., Griebel, M.: A representer theorem for deep kernel
  learning. Journal of Machine Learning Research  \textbf{20},  1--32 (2019)

\bibitem{grepl2005aposteriori}
Grepl, M., Patera, A.: A posteriori error bounds for reduced-basis
  approximations of parametrized parabolic partial differential equations.
  {ESAIM}: Mathematical Modelling and Numerical Analysis  \textbf{39}(1),
  157--181 (2005)

\bibitem{haasdonk2017reduced}
Haasdonk, B.: Reduced basis methods for parametrized {PDE}s -- a tutorial
  introduction for stationary and instationary problems. In: Model Reduction
  and Approximation: Theory and Algorithms, pp. 65--136. SIAM, Philadelphia
  (2017)

\bibitem{haasdonk2022certified}
Haasdonk, B., Kleikamp, H., Ohlberger, M., Schindler, F., Wenzel, T.: A new
  certified hierarchical and adaptive {RB}-{ML}-{ROM} surrogate model for
  parametrized {PDE}s. arXiv:2204.13454  (2022), to appear in SISC 2023.

\bibitem{himpe2018hierarchical}
Himpe, C., Leibner, T., Rave, S.: Hierarchical approximate proper orthogonal
  decomposition. SIAM J. Sci. Comput.  \textbf{40}(5),  A3267--A3292 (2018)

\bibitem{santin2021kernel}
Santin, G., Haasdonk, B.: Kernel methods for surrogate modeling. In: Model
  Order Reduction, vol.~2. De Gruyter (2021)

\bibitem{wang2019nonintrusive}
Wang, Q., Hesthaven, J.S., Ray, D.: Non-intrusive reduced order modeling of
  unsteady flows using artificial neural networks with application to a
  combustion problem. Journal of Computational Physics  \textbf{384},  289--307
  (2019)

\bibitem{wenzel2023data}
Wenzel, T., Marchetti, F., Perracchione, E.: Data-driven kernel designs for
  optimized greedy schemes: A machine learning perspective. arXiv:2301.08047
  (2023)

\bibitem{wenzel2021universality}
Wenzel, T., Santin, G., Haasdonk, B.: Universality and {O}ptimality of
  {S}tructured {D}eep {K}ernel {N}etworks. arXiv:2105.07228  (2021)

\end{thebibliography}
\end{document}